\documentclass{amsart}
\usepackage[utf8]{inputenc}

\usepackage{graphics}
\usepackage{thmtools}
\usepackage{wasysym}
\usepackage[T1]{fontenc}    

\usepackage{amsthm}
\usepackage{amsbsy,amsmath,amssymb,amscd,amsfonts}
\usepackage[pagebackref=true]{hyperref}
\usepackage[nameinlink,capitalize,noabbrev]{cleveref}

\usepackage{graphicx,float,latexsym,color}
\usepackage[export]{adjustbox}
\usepackage[font={scriptsize,it}]{caption}
\usepackage{subcaption}

\usepackage{makecell}

\usepackage[dvipsnames]{xcolor}

\newtheorem*{theorem*}{Theorem}
\newtheorem{observation}{Observation}
\newtheorem{proposition}{Proposition}

\newtheorem{corollary}{Corollary}
\newtheorem{lemma}{Lemma}
\theoremstyle{remark}

\theoremstyle{definition}
\newtheorem{definition}{Definition}

\hypersetup{
    pdftoolbar=true,        
    pdfmenubar=true,        
    pdffitwindow=false,     
    pdfstartview={FitH},    
    colorlinks=true,       
    linkcolor=OliveGreen,          
    citecolor=blue,        
    filecolor=black,      
    urlcolor=red           
}

\usepackage{listings}
\usepackage{lineno}

\newcommand{\ti}[1]{\textit{#1}}

\usepackage{comment}
\usepackage{microtype}
\usepackage{footnote}

\renewcommand{\O}{\mathcal{O}}
\newcommand{\B}{\mathcal{B}}
\newcommand{\E}{\mathcal{E}}

\newcommand{\T}{\mathcal{T}}

\renewcommand{\L}{\mathcal{L}}

\renewcommand{\T}{\mathcal{T}}

\newcommand{\ol}{\overline}
\renewcommand{\l}{\lambda}

\title[Stationary Kiepert focus]{The stationary focus of the Kiepert parabola over a special Poncelet triangle family}
\author[M. Helman]{Mark Helman}
\author[R. Garcia]{Ronaldo A. Garcia}
\author[D. Reznik]{Dan Reznik} 

\begin{document}

\begin{abstract}
We show that the focus of the Kiepert in-parabola remains stationary over a family of circle-inscribed Poncelet triangles which contain an equilateral triangle. 
\end{abstract}

\maketitle

\section{Introduction}
\label{sec:intro}
The Kiepert in-parabola is tangent to the three sides of a triangle \cite{mw}. We show that its focus -- called $X_{110}$ in \cite{etc} -- remains stationary over a family of circle-inscribed Poncelet triangles which contain an equilateral triangle.

In a companion article we study degeneracies of a closely-related family, namely, Poncelet triangles circumscribing their incircle, also containing one equilateral member  \cite[Sec.4]{garcia2024-incircle}.

\begin{figure}
\centering
\includegraphics[width=0.7\linewidth]{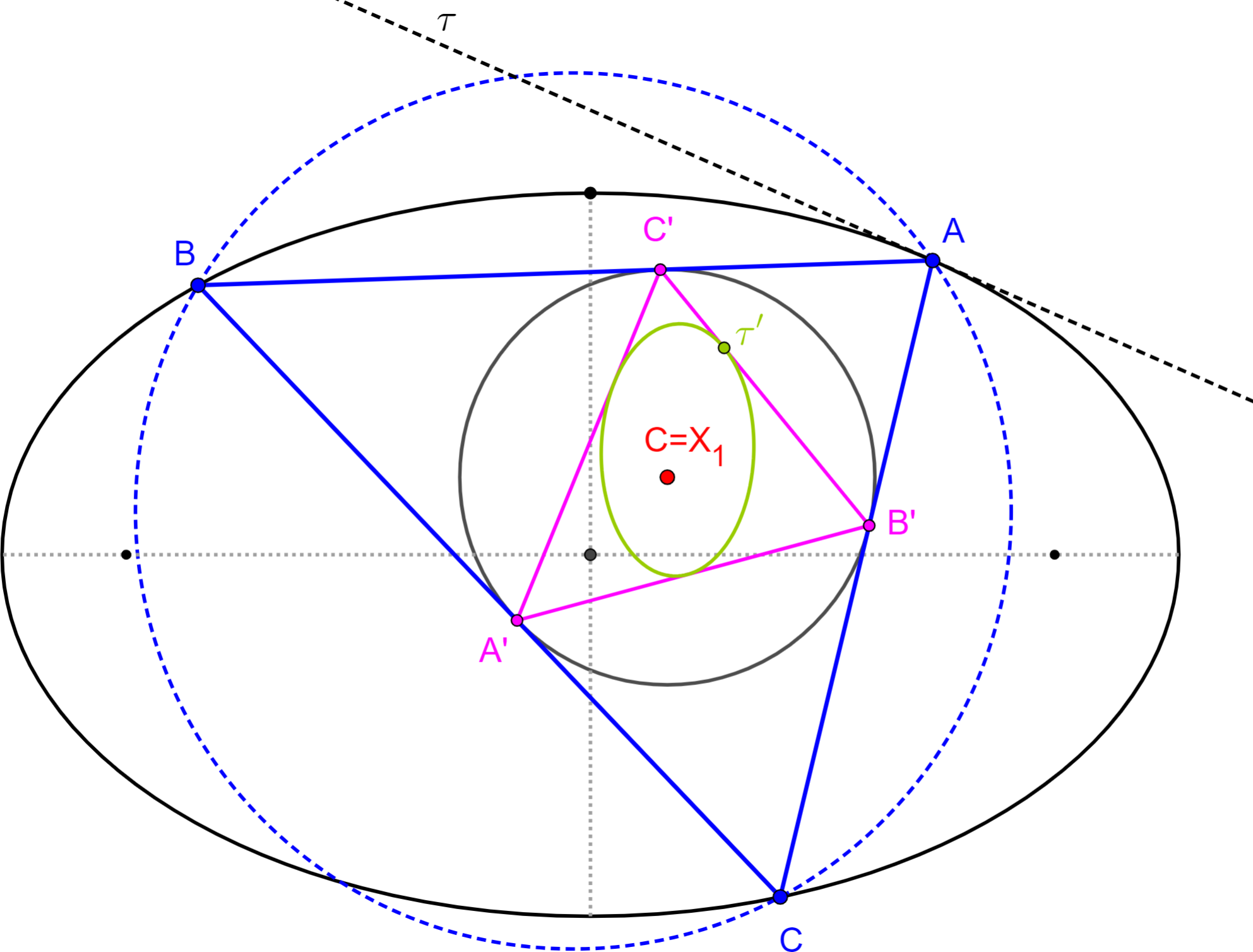}
\caption{The contact family  $A'B'C'$ is circle-inscribed and envelops a conic which is the polar image of tangents $\tau$ to $\E$ with respect to the incircle. The point $\tau'$ indicated the instantaneous polarity.}
\label{fig:contact-dual}
\end{figure}

As shown in \cref{fig:contact-dual}, the two families mentioned are polar images of each other with respect to their outer conics \cite{akopyan2007-conics}.

A key behavior described here involves the stationarity of the focus $X_{110}$ (using \cite{etc} terminology) of the \ti{Kiepert Parabola}, an well-studied inconic to a triangle, see \cite{mw}. This corresponds to the stationarity of Feuerbach's point of Poncelet triangles about the incircle, described in \cite{garcia2024-incircle}, but also proved here. 

\subsection{Computational proofs}

Most phenomena visited have been discovered via simulation. In some figure captions we include links to videos as well as interactive animations with a custom-built tool  \cite{darlan2021-app}.

\section{Review: Poncelet's porism}
\label{sec:poncele}
Poncelet's porism is a 1d family of $n$-gons with vertices on a first conic $\E$ and with sides tangent to a second conic $\E_c$ (also called the `caustic'). For such a porism to exist, $\E$ and $\E_c$ must be positioned in $\mathbb{R}^2$ so as to satisfy `Cayley's condition' \cite{dragovic11}. While the porism is projectively invariant -- a conic pair is the projective image of two circles where the `Poncelet map' is linearized, -- we have found that Poncelet porisms of triangles ($n=3$) are a wellspring of interesting Euclidean phenomena involving the dynamical geometry of classical objects associated with the triangle (centers, circles, lines/axes, etc.).

An emblematic case is when when $\E,\E_c$ are \ti{confocal}. In such a porism, -- also known as the \ti{elliptic billiard} -- the perimeter is conserved \cite{sergei91}, and the loci of all four `Greek centers' of the triangle -- incenter, barycenter, circumcenter, orthocenter -- sweep ellipses \cite{fierobe2021-circumcenter,garcia2020-ellipses,olga14}. Furthermore, the family conserves the ratio of inradius-to-circumradius (equivalently, the sum of internal angle cosines) \cite{garcia2020-new-properties}, a result which was extended to all $n$ \cite{akopyan2020-invariants}. 


\section{Equilateral containment}
\label{sec:equilateral}
We study the case of circular $\E$ and $\E_c$ is a nested ellipse. This is the affine image of a Poncelet triangle family inscribed in some conic and circumscribing a circular $\E_c$, and can therefore be described via the symmetric parametrization of \cite[Sec.2]{garcia2024-incircle}. In particular, we are interested in the circle-inscribed family which contains an equilateral member. 

\begin{observation}
A family of Poncelet triangles inscribed in a circle $\E$ contains an equilateral triangle if and only if its polar image with respect to $\E$ also contains one.
\end{observation}

\begin{figure}
\centering
\includegraphics[width=0.8\linewidth]{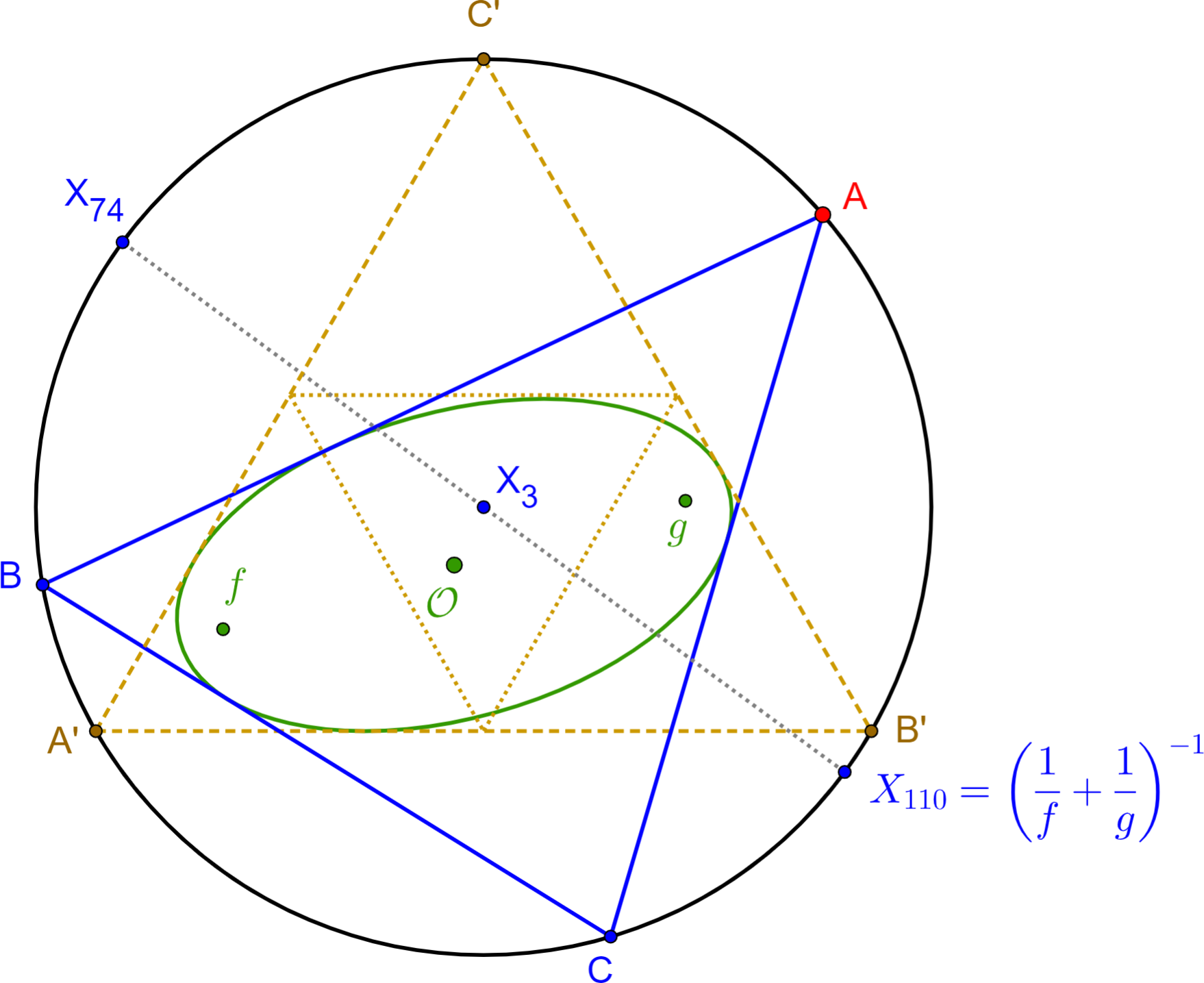}
\caption{$T_o'=A'B'C'$ (brown) is an equilateral for which an inconic (green), centered at $\O$, is chosen, with foci at (complex) $f,g$. $ABC$ is a triangle (blue) in the Poncelet family defined by the circumcircle (black) and the chosen inconic. Over the family, $X_{110}$. the focus of the Kiepert parabola (not shown) is stationary at $(1/f+1/g)^{-1}$, as is its antipode $X_{74}$, i.e., called in \cite{etc} the `isogonal conjugate of the Euler infinity point'. \href{https://youtu.be/czK2ZQycq24}{Video}}
\label{fig:ci-basic}
\end{figure}

Let $\E_c$ be an inconic with foci $f,g\in\mathbb{C}$, such that there is a Poncelet family of triangles inscribed in the unit circle $\mathbb{T}$ and circumscribing $\E_c$. Referring to \cref{fig:ci-basic}:

\begin{lemma}
A Poncelet family of triangles inscribed in the unit circle and circumscribing an ellipse with foci $f,g\in\mathbb{C}$ contains an equilateral triangle iff $1/f+1/g$ is on the unit circle $\mathbb{T}$.
\label{lem:equi}
\end{lemma}

\begin{proof}
First, assume that $T_o$ is an equilateral triangle in this family. Let $\zeta\in\mathbb{T}$ be one of the vertices of $T_o$ such that the other vertices are $e^{\frac{2\pi}{3}i}\zeta$ and $e^{\frac{4\pi}{3}i}\zeta$. Using the symmetric parameterization described in \cite[Sec.2]{garcia2024-incircle}, let $\lambda_o\in\mathbb{T}$ be the parameter corresponding to the triangle $T_o$. Plugging in the vertices in the parameterization, we get
\begin{gather*}
    f+g+\l_o\ol f \ol g = \zeta+e^{\frac{2\pi}{3}i}\zeta+e^{\frac{4\pi}{3}i}\zeta= \zeta\left(1+e^{\frac{2\pi}{3}i}+e^{\frac{4\pi}{3}i}\right)=0\\
    \implies \l_o=-\frac{(f+g)}{\ol f\ol g}\implies \left|-{\frac{(f+g)}{\ol f\ol g}}\right|=|\l_o|=1\\
    \implies |f+g|=|f g|\implies\left|\frac1f+\frac1g\right|=\left|\frac{f+g}{f g}\right|=1
\end{gather*}
so $1/f+1/g$ is on the unit circle $\mathbb{T}$, as desired.

Conversely, suppose that $1/f+1/g$ is on the unit circle, that is, $|1/f+1/g|=1$. Rearranging, we have $|f+g|=|f g|$. Define $\l_o:=-(f+g)/(\ol f\ol g)$. Then
\[
|\l_o|=\left|-\frac{(f+g)}{\ol f\ol g}\right|=\frac{|f+g|}{|\ol f\ol g}=\frac{|f+g|}{|f g|}=1
\]

By definition, $f+g+\l_o\ol f \ol g =0$. Moreover,
\begin{gather*}
    f g +\l_o (\ol f+\ol g)=\l_o \ol{\l_o}f g +\l_o (\ol f+\ol g)=\l_o (\ol{\l_o} f g+\ol f+\ol g)=\\
    \l_o \ol{\left(f+g+\l_o \ol f\ol g\right)}=\l_o \ol 0=0
\end{gather*}

Hence, letting $z_1,z_2,z_3$ take the values of the 3 roots of $\zeta^3=\l_o$, we see that these $3$ roots form an equilateral triangle that satisfies the equations of the symmetric parameterization of \cite[Sec.2]{garcia2024-incircle}. Thus, this equilateral triangle is part of the Poncelet family, as desired.
\end{proof}

\section{Stationary focus of the Kiepert Parabola}
\label{sec:kipert}
\begin{figure}
\centering
\includegraphics[width=0.6\linewidth]{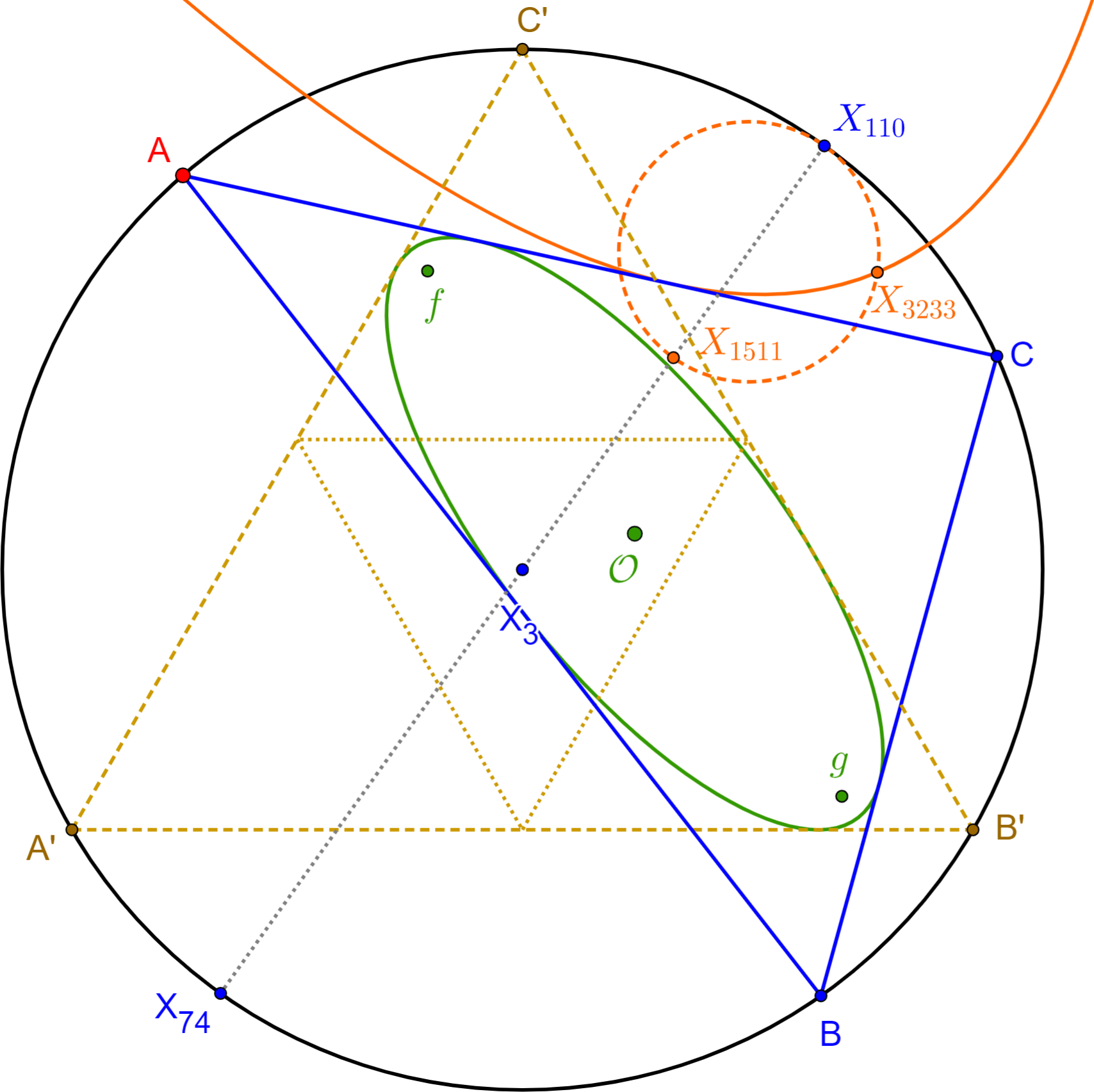}
\caption{A $T_o=ABC$ in $\T_o$ (blue) is shown (blue) along with its Kiepert parabola (orange), also an inconic, with focus at $X_{110}$. Over $\T_o$, the locus of its vertex $X_{3233}$ is a circle (dashed orange) whose diameter is $X_{110} X_{1511}$. \href{https://youtu.be/5QK-JzN16tM}{Video}}
\label{fig:ci-kiepert-vtx}
\end{figure}

Referring to \cref{fig:ci-kiepert-vtx}:

\begin{definition}
The \ti{Kiepert parabola} of a triangle is tangent to the three sides of a triangle. Its focus is $X_{110}$, always on the circumcircle, and its directrix is the Euler line $X_2 X_3$ \cite{mw}.
\end{definition}

Henceforth, let $\T_o^*$ denote a Poncelet family of triangles interscribed between the unit circle and some inconic $\E_c$ satisfying \cref{lem:equi}. Still referring to \cref{fig:ci-basic}, it follows from \cref{lem:equi}:

\begin{proposition}
Over $\T_o^*$, $X_{110}$ is stationary at $(1/f+1/g)^{-1}$.
\end{proposition}

On \cite{etc}, (i) the midpoint of $X_3$ and $X_{110}$ is called $X_{1511}$, and (ii) $X_{3233}$ is the vertex of the Kiepert inparabola \cite{moses2024-private}. Referring to \cref{fig:ci-kiepert-vtx}:

\begin{observation}
Over $\T_o^*$, $X_{3233}$ sweeps a circle with diameter $X_{110} X_{1511}$, i.e., $1/2$.  
\end{observation}

\begin{figure}
\centering
\includegraphics[width=0.8\linewidth]{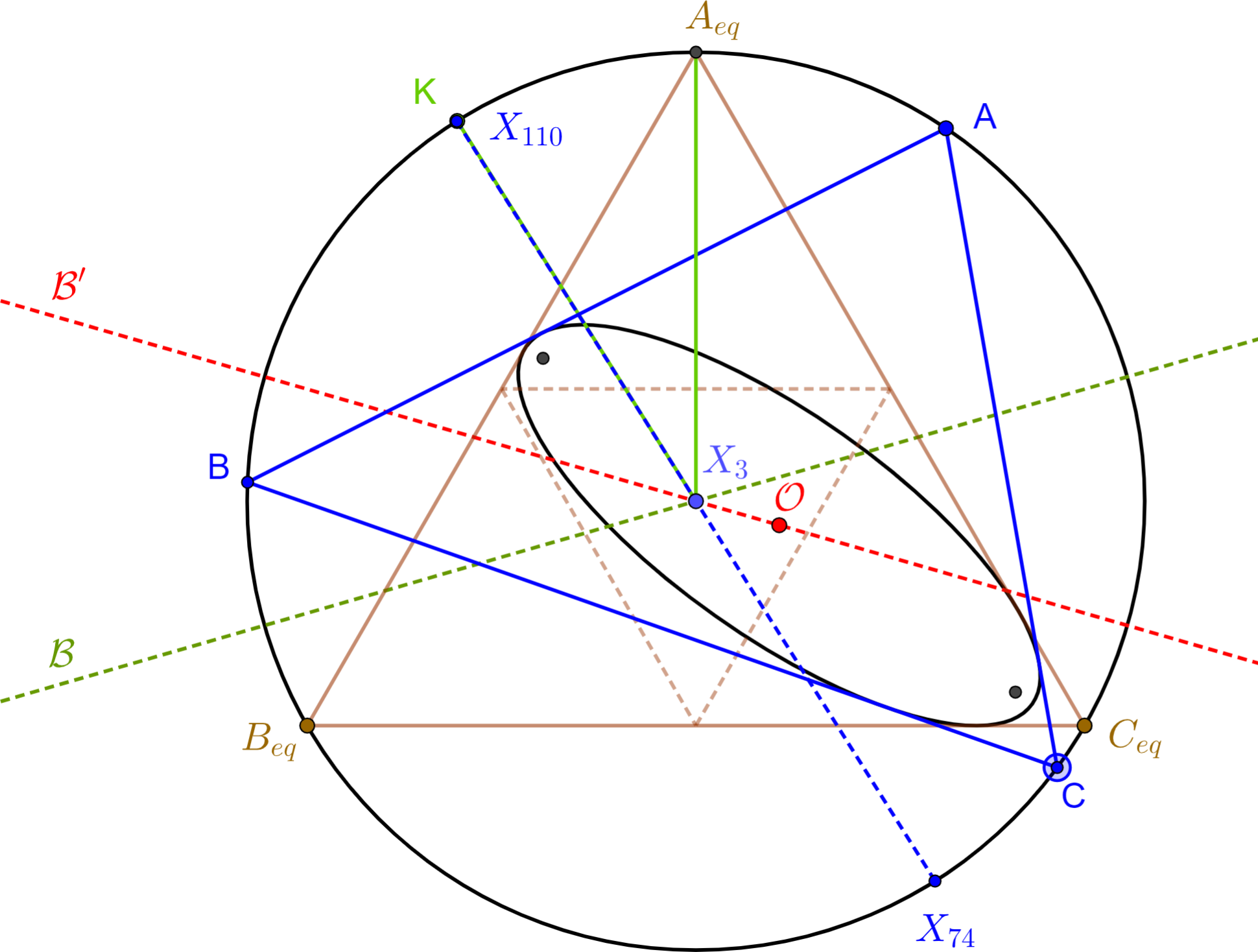}
\caption{The construction in \cref{prop:double-inv-1}: $K$ is chosen on the circumcircle; $A_{eq}$ is a vertex of the equilateral; $T_o=ABC$ is a Poncelet triangle in $\T_o^*$; $\B$ (dashed green) is the external bisector of $\angle K X_3 A_{eq}$; $\B'$ (dashed red) is its reflection about $X_3 A_{eq}$; over all caustic centers $\O$ on $\B'$, $X_{110}$ is stationary, over $\T_o^*$, at the chosen $K$.}
\label{fig:double-x110}
\end{figure}

Referring to \cref{fig:double-x110}, 
let $A_{eq}$ be one of the vertices of the equilateral in $\T_o^*$. Let $K$ be a chosen point on the circumcircle. Let $\B$ represent the external bisector of $\angle K X_3 A_{eq}$. It can be shown:

\begin{proposition}
\label{prop:double-inv-1}
Over $\T_o^*$, $X_{110}$ will be stationary at $K$ for any $\O$ on the reflection of $\B$ about $X_3 A_{eq}$.
\end{proposition}

\begin{figure}
\centering
\includegraphics[width=0.8\linewidth]{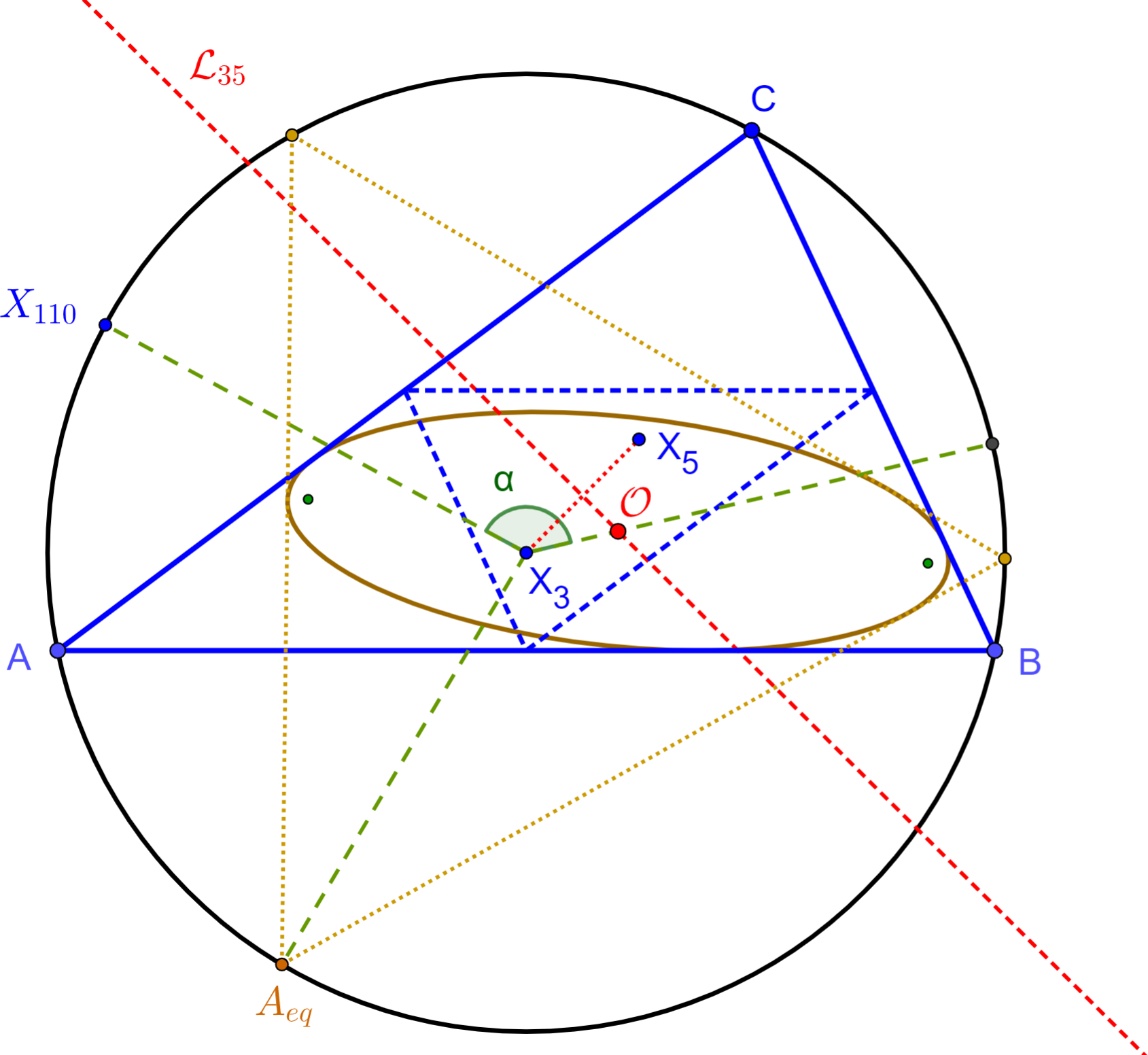}
\caption{The construction of \cref{prop:l35,prop:l35-aeq}. $X_{110}$ (of $T_o=ABC$) is stationary over Poncelet if $\O$ lies on $\L_{35}$ (dashed red), the perpendicular bisector of $X_3$ and $X_5$. A vertex $A_{eq}$ of the equilateral contained in the family is shown.}
\label{fig:l35}
\end{figure}

Referring to \cref{fig:l35}, let $T_o=ABC$ be some generic triangle. Let $\O$ be the center of an inconic $\E_c$. Let $\T_o$ be the Poncelet family defined by the circumcircle and $\E_c$. Let $\L_{35}$ be the perpendicular bisector of $X_3 X_5$. It can be shown:

\begin{proposition}
\label{prop:l35}
$\T_o$ will contain an equilateral iff $\O$ lies on $\L_{35}$. Furthermore, the stationary $X_{110}$ is independent of $\O$. 
\end{proposition}

Let $\alpha=\angle X_{110} X_3 \O$. Let $\O'$ be the counterclockwise rotation of $\O$ about $X_3$ by $\alpha/3+\pi$. Still referring to \cref{fig:l35}, it can be shown:

\begin{proposition}
\label{prop:l35-aeq}
One vertex of the equilateral contained in $\T_o$ is located at the intersection of $X_3 \O'$ with the circumcircle.
\end{proposition}

\section{Stationary Feuerbach point of the polar}
\label{sec:feuerbach}

Let $\T_o'$ be the family of tangential triangles to $\T_o$, i.e., whose sides are tangent to the circumcircle at the vertices of $\T_o$ \cite[Tangential Triangle]{mw}. Clearly, this family circumscribes a circle and will be inscribed in an ellipse $\E$, called above $\T$.

\begin{observation}
$\E$ is the polar image of the chosen inconic with respect to the circumcircle.
\end{observation}

\begin{corollary}
Since the tangential of an equilateral is an equilateral, $\T_o'$ will contain an equilateral, and therefore the center of its incircle will be on $\E_{eq}$.
\end{corollary}

In turn, this leads to a restatement of \cite[Prop.15]{garcia2024-incircle}: 

\begin{corollary}
The Feuerbach point $X_{11}$ of $\T_o'$ will be stationary.
\end{corollary}

\begin{proof}
$X_{110}$ of a reference triangle is $X_{11}$ of its tangential. Conversely, $X_{11}$ of the reference is $X_{110}$ of the contact triangle.
\end{proof}

On \cite{etc}, $X_{65}$ is the orthocenter $X_4$ of the contact triangle.

\begin{corollary}
The locus of $X_{65}$ of the tangential family is a circle.  
\end{corollary}

\begin{proof}
The contact family is Poncelet (sides envelop the polar image of $\E$ wrt incircle) and circle-inscribed. $X_4$ is always homothetic to a 90-degree-rotated version of the inscribing conic. 
\end{proof}

\bibliographystyle{maa}
\bibliography{refs_01_pub,refs_03_sub,refs_04_unsub,refs}

\end{document}